\documentclass[english]{smfart}
\usepackage{amsmath,amsfonts,amsthm,amssymb,amscd,color,xcolor,mathrsfs}
\usepackage{hyperref}

\author[Armen Shirikyan]{Armen Shirikyan}
\address{University of Cergy-Pontoise (France) and 
National Research University MPEI, Moscow (Russia)}

\email{\href{mailto:Armen.Shirikyan@u-cergy.fr}{Armen.Shirikyan@u-cergy.fr}}
\title[Mixing for the stochastic Burgers equation]{Mixing for the Burgers equation driven by a localised two-dimensional stochastic forcing}

\colorlet{darkblue}{blue!50!black}

\hypersetup{
    colorlinks,%
    citecolor=darkblue,%
    filecolor=red,%
    linkcolor=darkblue,%
    urlcolor=magenta,%
    pdfnewwindow=true,%
    pdfstartview={FitH}
}

\newcommand{\dd}{{\textup d}}

\newcommand{\p}{\partial}
\newcommand{\e}{\varepsilon}

\newcommand{\E}{{\mathbb E}}
\newcommand{\R}{{\mathbb R}}
\newcommand{\IP}{{\mathbb P}}

\newcommand{\BBB}{{\boldsymbol B}}

\newcommand{\oomega}{{\boldsymbol\omega}}
\newcommand{\ttau}{{\boldsymbol\tau}}

\newcommand{\OOmega}{{\boldsymbol\Omega}}

\newcommand{\BB}{{\mathcal B}}

\newcommand{\LL}{{\mathcal L}}

\newcommand{\PP}{{\mathcal P}}

\newcommand{\XX}{{\mathcal X}}
\newcommand{\YY}{{\mathcal Y}}

\newcommand{\PPPP}{{\mathfrak P}}

\newcommand{\uuu}{{\boldsymbol{\mathit u}}}
\newcommand{\vvv}{{\boldsymbol{\mathit v}}}

\newcommand{\VVV}{{\boldsymbol{\mathit V}}}

\newcommand{\lspan}{\mathop{\rm span}\nolimits}

\theoremstyle{plain}

\newtheorem{theorem}{Theorem}[section]
\newtheorem{lemma}[theorem]{Lemma}
\newtheorem{proposition}[theorem]{Proposition}

\theoremstyle{definition}

\theoremstyle{remark}

\begin{document}
\frontmatter
\begin{abstract}
We consider the one-dimensional Burgers equation perturbed by a stochastic forcing, which is assumed to be white in time and localised and low-dimensional in space. We establish a mixing property for the Markov process associated with the problem in question. The proof is based on a general criterion for mixing and a recent result on global approximate controllability to trajectories for damped conservation laws. 
\end{abstract}


\thanks{This research was supported by the RSF grant 14-49-00079}
\maketitle

\mainmatter
%\tableofcontents

\section{Introduction}
The paper is devoted to studying the problem of uniqueness and mixing of a stationary measure for the following stochastic Burgers equation:
\begin{align}
\p_tu-\nu\p_x^2u+u\p_xu&=h(x)+\sum_{j=1}^2\dot \beta_j(t)e_j(x), \quad x\in I,\label{0.1}\\
u(t,0)=u(t,\pi)&=0, \label{0.2}
\end{align}
where $I=(0,\pi)$, $\nu>0$ is a parameter, $h\in L^\infty(I)$ is a fixed function, $\{\beta_j\}$ are independent Brownian motions, and~$\{e_j\}$ are continuous functions supported by an interval $[a,b]\subset I$. 
The Cauchy problem for~\eqref{0.1}, \eqref{0.2} is well posed (see Chapter~14 in~\cite{DZ1996}): for any $u_0\in L^2(I)$ there is a unique stochastic process $u(t,x)$ whose almost every trajectory coincides with~$u_0$ for $t=0$ and satisfies the inclusion
\begin{equation} \label{0.3}
u\in \XX:=C(\R_+,L^2(I))\cap L_{\mathrm{loc}}^2(\R_+,H_0^1(I)) 
\end{equation}
and Eq.~\eqref{0.1} in the sense of distributions; see the notation section below for the definition of functional spaces. Moreover, the family of all solutions form a  Markov process in~$L^2(I)$, which possesses the Feller property. The latter holds also in the more regular space~$H^1_0(I)$.  The Bogolyubov--Krylov argument enables one to prove the existence of a stationary measure; see Chapter~14 in~\cite{DZ1996}. Our aim is to study the uniqueness and mixing of a stationary measure. 

The main result of this article proves that, for a particular choice of the functions~$e_1$ and~$e_2$, the Markov process associated with problem~\eqref{0.1}, \eqref{0.2} has a unique stationary measure~$\mu_\nu$ for any $\nu>0$, and the law of any solution converges weakly to~$\mu_\nu$ as $t\to+\infty$ in the Kantorovich--Wasserstein metric. An exact formulation of this result is given in Section~\ref{s2.1}, and the proof is presented in Section~\ref{s3}. 

Let us mention that the problem of mixing of dissipative PDE's with a random external force was in the focus of attention of many researchers in the last twenty years, and it is rather well understood in the case when when all deterministic modes are forced; see the literature review in Chapter~3 of~\cite{KS-book}. The situation in which the support of the random forcing does not cover all deterministic modes is far less studied, and only few results are available. Namely, Hairer and Mattingly~\cite{HM-2006,HM-2011} established the uniqueness of a stationary measure and its exponential stability for 2D Navier--Stokes equations with a random forcing, white in time and finite-dimensional in space, provided that the equation is considered on a torus or a sphere, and the unperturbed problem has a unique globally stable equilibrium. Bortichev~\cite{boritchev-2013} considered the viscous Burgers equation on a circle and proved that, if the unperturbed dynamics is globally asymptotically stable, then for any random perturbation, the stationary measure is unique, and convergence to it holds with a polynomial rate independent of the viscosity; see also~\cite{EKMS-2000,IK-2003,DV-2015} for the case of an inviscid equation. The paper~\cite{shirikyan-aens2015}  establishes the property of exponential mixing for the 2D Navier--Stokes system in a bounded domain in the case when the noise is smooth, bounded, and localised in the physical space and in time. Finally, F\"oldes et\,al~\cite{FGRT-2015} proved that the Hairer--Mattingly result is true for the Boussinesq system perturbed by a highly-degenerate stochastic forcing which acts only on the equation for temperature. To the best of our knowledge, the situation in which the unperturbed dynamics is non-trivial, while the noise is low-dimensional and localised in space was not considered earlier, and this is the main novelty of this work. At the same time, it remains an open question to determine the rate of convergence to the stationary measure. 

\subsection*{Notation}
For a separable Banach space~$X$, we denote by $\BB_X$ the Borel $\sigma$-algebra on~$X$, by~$C_b(X)$ the space of bounded continuous functions $f:X\to \R$ endowed with the supremum norm~$\|\cdot\|_\infty$, and by $L_b(X)$ the subspace of functions $f\in C_b(X)$ such that
$$
\|f\|_L:=\|f\|_\infty+\sup_{u\ne v}\frac{f(u)-f(v)|}{\|u-v\|_X}<\infty. 
$$
The set of Borel probability measures on~$X$ is denoted by~$\PP(X)$ and endowed with the Kantorovich--Wasserstein distance (which metrises the weak convergence of measures) defined as
$$
\|\mu-\nu\|_L^*=\sup_{\|f\|_L\le1}\bigl|\langle f,\mu\rangle-\langle f,\nu\rangle\bigr|,
$$
where $\langle f,\mu\rangle$ stands for the integral of~$f$ over~$X$ with respect to~$\mu$. Given $R>0$, we write $B_X(R)$ for the closed ball in~$X$ centred at zero of radius~$R$. 

Let $I=(0,\pi)$ and let $J\subset\R$ be a closed interval. We use the following functional spaces:

$C(J,X)$ is the space of continuous functions $f:J\to X$.

$L^2(J,X)$ is the space of Borel-measurable functions $f:J\to X$ such that $\|f(t)\|_X$ is square-integrable on~$J$. If $J$ is unbounded, then $L_{\mathrm{loc}}^2(J,X)$ stands for the space of $X$-valued functions on~$J$ whose restriction to any finite interval $J'\subset J$ belongs to~$L^2(J',X)$.

$H^s(I)$ is the Sobolev space of order~$s\ge0$ with a standard norm~$\|\cdot\|_s$ and $H_0^s(I)$ is the closure in~$H^s(I)$ of the space of infinitely smooth functions with compact support. We write $H=L^2(I)$, $V=H_0^1(I)$, $U=H^2(I)\cap V$, and denote by~$\|u\|$ and~$\|u\|_V=\|\p_xu\|$ the norms in~$H$ and~$V$, respectively.

$\LL(X,Y)$ stands for the space of continuous linear operators from~$X$ to~$Y$. 

We denote by~$C_i$, $i=1,2,\dots$, unessential positive numbers. 

\section{Main result}
\subsection{The result}
\label{s2.1}
We consider problem~\eqref{0.1}, \eqref{0.2} on the interval~$I=(0,\pi)$. Under the regularity hypotheses on~$e_1$, $e_2$, and~$h$ mentioned in the Introduction, this problem is globally well posed and generates a Feller family of Markov processes, which will be denoted by~$(u_t,\IP_u)$. Let $P_t(u,\Gamma)=\IP_u\{u_t\in\Gamma\}$ be its transition function and let
$$
\PPPP_t:C_b(H)\to C_b(H), \quad \PPPP_t^*:\PP(H)\to\PP(H)
$$
be the corresponding Markov semigroups. Recall that a measure~$\mu\in\PP(H)$ is said to be {\it stationary\/} for $(u_t,\IP_u)$  if $\PPPP_t^*\mu=\mu$ for $t\ge0$. 

Let us fix an interval $[a,b]\subset I$ and assume that the functions~$e_j$ entering~\eqref{0.1} have form
$$
e_1(x)=c_1\sin\Bigl(\pi\frac{x-a}{b-a}\Bigr), \quad
e_2(x)=c_2\sin\Bigl(2\pi\frac{x-a}{b-a}\Bigr) \quad\mbox{for $x\in[a,b]$}, 
$$
where $c_1,c_2\in\R$ are non-zero numbers, and $e_1(x)=e_2(x)=0$ for $x\notin[a,b]$. The following theorem is the main result of this paper. 

\begin{theorem} \label{t1}
For any $\nu>0$, the Markov family~$(u_t,\IP_u)$ associated with~\eqref{0.1}, \eqref{0.2} has a unique stationary measure~$\mu_\nu\in\PP(H)$. Moreover, the measure~$\mu_\nu$ is concentrated on~$V$, and there is a function $\alpha(t)$, defined on the positive half-line and going to zero as $t\to+\infty$, such that 
\begin{equation} \label{4}
\|\PPPP_t^*\lambda-\mu_\nu\|_L^*\le\alpha(t)\quad\mbox{for $t\ge1$}, 
\end{equation}
where $\lambda\in\PP(H)$ is an arbitrary initial measure, and $\|\cdot\|_L^*$ is the Kantorovich--Wasserstein distance considered over the space~$V$. 
\end{theorem}

Let us emphasise that the convergence in~\eqref{4} (called {\it mixing property\/} for~$\mu_\nu$) holds uniformly with respect the initial measures~$\lambda$. This is due to the strong dissipative character of the Burgers equation; cf.~\cite{boritchev-2013}. 

A proof of Theorem~\ref{t1} is given in the next section. Here we present the general idea of the proof and outline the main steps. 

\subsection{Scheme of the proof of Theorem~\ref{t1}}
\label{s2.2}
As was mentioned in the Introduction, the existence of a stationary measure~$\mu_\nu$ follows from the Bogolyubov--Krylov argument and standard a priori estimates for solutions of the Burgers equation; e.g., see Section~2.5 in~\cite{KS-book} for the more complicated case of the Navier--Stokes system. Thus, we shall concentrate on the proof of uniqueness of~$\mu_\nu$ and mixing property.

\smallskip
{\it Step~1. Reduction to regular solutions.\/} 
The Burgers equation possesses the following regularising property: 
\begin{description}
\item[Regularisation] 
For any $\lambda\in\PP(H)$ and $t>0$, we have $(\PPPP_t^*\lambda)(V)=1$. 
\end{description}
It follows that any stationary measure is concentrated on~$V$, and it suffices to prove uniqueness of stationary distribution in the class of probability measures  on~$V$. Furthermore, when proving convergence~\eqref{4}, one can assume that $\lambda(V)=1$. Finally, recalling the relation
$$
(\PPPP_t^*\lambda)(\Gamma)=\int_VP_t(v,\Gamma)\lambda(\dd v), \quad \Gamma\in\BB_H,
$$
we see that it suffices to prove~\eqref{4} with $\lambda=\delta_v$ with an arbitrary $v\in V$:
\begin{equation} \label{5}
\|P_t(v,\cdot)-\mu_\nu\|_{L}^*\le\alpha(t)\quad\mbox{for $t\ge 1$, $v\in V$}. 
\end{equation}

\smallskip
{\it Step~2. A sufficient condition for uniform mixing\/}. 
We shall consider the restriction of the Markov process~$(u_t,\IP_u)$ to the space~$V$. Let us define the Markov process $(\uuu_t,\IP_\uuu)$ in~$\VVV=V\times V$ as a pair of independent copies of~$(u_t,\IP_u)$. In other words, the probability space is defined as the direct product $\OOmega=\Omega\times\Omega$, with the natural product $\sigma$-algebra on it, the process~$\uuu_t$ is given by
$$
\uuu_t(\oomega)=\bigl(u_t(\omega),u_t(\omega')\bigr), \quad 
\oomega=(\omega,\omega')\in\OOmega,
$$
and $\IP_\uuu=\IP_u\otimes\IP_{u'}$ for $\uuu=(u,u')$. 
In view of Theorem~3.1.3 in~\cite{KS-book} and the remark following it, the uniqueness of a stationary measure and convergence~\eqref{5} will be established if we prove the following two properties:
\begin{description}
\item[Stability]
There is a sequence $\{B_m\}$ of closed subsets of~$V$ such that
\begin{equation} \label{6}
\sup_{t\ge1}\|P_t(v,\cdot)-P_t(v',\cdot)\|_{L}^*\le\delta_m\quad 
\mbox{for $v,v'\in B_m$},
\end{equation}
where $\{\delta_m\}$ is a sequence of positive numbers going to zero as $m\to\infty$, and~$\|\cdot\|_L^*$ stands for the Kantorovich--Wasserstein distance over~$V$. 

\smallskip
\item[Recurrence]
Let $\ttau_m$ be the first hitting time of the set $\BBB_m=B_m\times B_m$ for the process $(\uuu_t,\IP_\uuu)$. Then 
\begin{equation} \label{7}
\sup_{\uuu\in\VVV}\IP_\uuu\{\ttau_m>t\} \le p(m,t)
\quad\mbox{for any $m\ge1$, $t>0$},
\end{equation}
where $p(m,t)\to0$ as $t\to+\infty$. 
\end{description}

The proofs of these assertions (outlined below) are based on the contraction of the $L^1$-norm of the difference between two solutions of the Burgers equation, a strong dissipation property and parabolic regularisation, and a recent controllability result for damped--driven conservation laws. 

\smallskip
{\it Step~3. Stability\/}. It is a well-known fact that, for any initial functions $v,v'\in V$, the corresponding solutions of~\eqref{0.1}, \eqref{0.2} satisfy the inequality
\begin{equation} \label{8}
\|u(t)-u'(t)\|_{L^1(I)}\le \|v-v'\|_{L^1(I)}\quad\mbox{for $t\ge0$};
\end{equation}
see Section~3.3 in~\cite{hormander1997}.
On the other hand, the a priori estimates for solutions of~\eqref{0.1}, \eqref{0.2} imply that, if an initial function belongs to a bounded set $B\subset V$, then 
\begin{equation} \label{9}
\sup_{t\ge1}\E\,\|u(t)\|_{H^{s}}\le C_s(B),
\end{equation}
where $s\in(1,2)$ is arbitrary. 
Combining inequalities~\eqref{8} and~\eqref{9} with an interpolation inequality, we can easily prove~\eqref{6} with any sequence of bounded sets $B_m\subset V$ whose diameters in~$L^1$ go to zero. 

\smallskip
{\it Step~4. Recurrence\/}. 
A comparison theorem (see Section~2.2 in~\cite{AL-1983}) and regularisation for parabolic equations imply that the solutions of~\eqref{0.1}, \eqref{0.2} satisfy the inequality
\begin{equation} \label{10}
\IP_v\bigl\{\|u_2\|_{s}\le R\bigr\}\ge p\quad\mbox{for $v\in V$},
\end{equation}
where $s\in (1,2)$ is fixed, $p>0$ is a number not depending on~$v$, and $R>0$ is sufficiently large. On the other hand, using global controllability to trajectories of the Burgers equation (see Theorem~1.2 in~\cite{RS-2015}) and the continuous dependence of trajectories on the noise, one can prove that
\begin{equation} \label{11}
\IP_v\bigl\{u_{T_m}\in B_m\bigr\}\ge\e_m>0\quad
\mbox{for $v\in B_{V\cap H^s}(R)$},
\end{equation}
where $B_m\subset V$ are bounded closed sets whose $L^1$-diameters go to zero as $m\to\infty$, and~$T_m>0$ are some numbers. Combining~\eqref{10} and~\eqref{11} with the Markov property, we derive the required result. 

\section{Proof of the main result}
\label{s3}
As was explained in Section~\ref{s2.2}, Theorem~\ref{t1} will be established if we prove the regularisation, stability, and recurrence properties stated in Steps~1 and~2. This is done in the next three subsections. 

\subsection{Regularisation}
\label{s3.1}
We wish to prove that, if~$u(t,x)$ is the solution of~\eqref{0.1}, \eqref{0.2} issued from an initial point~$u_0(x)$, which is an $H$-valued random variables independent of~$(\beta_1,\beta_2)$, then
\begin{equation} \label{12}
\IP\{u(t)\in V\}=1\quad\mbox{for any $t>0$}. 
\end{equation}
We begin with a number of simple remarks. In view of the relation
$$
\IP\{u(t)\in V\}=\int_HP_t(z,V)\lambda(\dd z),
$$
where~$\lambda$ is the law of the initial state~$u_0$, it suffices to consider the case in which~$u_0$ is a deterministic function. Applying the Ito formula to the $L^2$ norm of a solution, using the $L^2$-orthogonality of the nonlinear term and solution, and repeating standard arguments (e.g., see the proof of Proposition~2.4.8 in~\cite{KS-book}), we derive
\begin{equation} \label{31}
\E\,\|u(t)\|^2+\E\int_0^t\|u(r)\|_V^2\,\dd r
\le C_1(1+t)\quad\mbox{for all $t\ge0$}.
\end{equation}
It follows that $\IP\{u(r)\in V\}=1$ for almost every $r\in(0,t)$. Thus, when proving~\eqref{12}, we can assume from the very beginning that~$u_0$ is a deterministic function belonging to~$V$.

\smallskip
Let us denote by~$z(t,x)$ the solution of the problem
\begin{equation} \label{13}
\p_tz-\nu\p_x^2z=\eta(t,x), \quad z(t,0)=z(t,\pi)=0, \quad z(0,x)=0,
\end{equation}
where~$\eta$ stands for the stochastic term in~\eqref{0.1}, and represent a solution of~\eqref{0.1}--\eqref{0.3} in the form $u=z+v$. In this case, $v(t,x)$ must be a solution of the problem
\begin{align}
\p_tv-\nu\p_x^2v+(v+z)\p_x(v+z)&=h(x), \label{14}\\
v(t,0)=v(t,\pi)&=0, \label{15}\\
v(0,x)&=u_0(x). \label{16}
\end{align}
The following lemma can be established by well-known arguments, and its proof is outlined at the end of this section. 

\begin{lemma} \label{l3.1}
Almost every trajectory of~$z$ belongs to the space
$$
\YY_s:=C(\R_+,V\cap H^s)\cap L_{\mathrm{loc}}^2(\R_+,U), \quad s<2,
$$ 
and there are positive numbers $C_k=C_k(\nu)$ such that
\begin{equation} \label{17}
\sup_{t\ge0}\E\Bigl(\|z(t)\|_1^2+\int_t^{t+1}\|z(r)\|_2^2\dd r\Bigr)^k\le C_k.
\end{equation}
Moreover, for any $\nu>0$, $k\ge1$, $T>0$, and $s\in(1,2)$, there is $C=C(k,T,s,\nu)>0$ such that
\begin{equation} \label{18}
\E\sup_{t\le r\le t+T}\|z(r)\|_s^k\le C\quad\mbox{for any $t\ge0$}. 
\end{equation}
\end{lemma}

In view of the continuous embedding $H^s\subset C(\bar I)$ valid for $s>1/2$ and inequality~\eqref{18}, for any $T>0$, the quantity 
\begin{equation} \label{28}
R_T:=\sup_{0\le t\le T}\bigl(\|z(t,\cdot)\|_{L^\infty}+\|\p_xz(t,\cdot)\|_{L^\infty}\bigr)
\end{equation}
is almost surely finite. Applying the maximum principle to~\eqref{14}--\eqref{16}, (e.g., see Section~2 in~\cite[Chapter~3]{landis1998}), we derive
\begin{equation} \label{26}
\|v\|_{L^\infty(D_T)}\le C_2(R_T+1)+\|u_0\|_{L^\infty},
\end{equation}
where we set $D_T=(0,T)\times I$. It follows that
\begin{equation} \label{27}
\|u\|_{L^\infty(D_T)}=\|v+z\|_{L^\infty(D_T)}\le (C_2+1)(R_T+1)+\|u_0\|_{L^\infty}. 
\end{equation}
Combining this with~\eqref{31}, we see that 
\begin{equation} \label{29}
g:=u\p_xu-h=(v+z)\p_x(v+z)-h\in L^2(D_T)\quad\mbox{almost surely}. 
\end{equation}
We now use the Duhamel formula to write a solution of~\eqref{14}-\eqref{16} in the form
\begin{equation} \label{30}
v(t,x)=(e^{tL}u_0)(x)-\int_0^t e^{(t-r)L}g(r,\cdot)\,\dd r,
\end{equation}
where $L=\nu \p_x^2$. Since the operator $e^{tL}$ is continuous from~$H$ to~$V$ for $t>0$, and its norm satisfies the inequality $\|e^{t L}\|_{\LL(H,V)}\le C_3t^{-1/2}$, we see that, with probability~$1$, $v(t)$ is a continuous function of $t>0$ with range in~$V$. Recalling that, by Lemma~\ref{l3.1}, the same is true for~$z(t)$, we arrive at~\eqref{12}. 

\subsection{Stability}
\label{s3.2}
We first show that inequalities~\eqref{8} and~\eqref{9} imply the stability property with any sequence of closed sets $B_m\subset V$ such that
\begin{equation} \label{33}
\sup_{v\in B_m}\|v\|_V\le C_4\quad\mbox{for any $m\ge1$},\qquad
\sup_{v,v'\in B_m}\|v-v'\|_{L^1}\to0\quad\mbox{as $m\to\infty$},
\end{equation}
where $C_4>0$ does not depend on~$m$. Indeed, given $v,v'\in B_m$, we denote by~$u(t)$ and~$u'(t)$ the solutions of~\eqref{0.1}--\eqref{0.3} with $u_0=v$ and $u_0=v'$, respectively. If $f\in L_b(V)$ is a function with $\|f\|_L\le1$, then 
\begin{equation} \label{34}
\bigl|\langle f,P_t(v,\cdot)\rangle -\langle f,P_t(v',\cdot)\rangle\bigr|
=\bigl|\E\,\bigl(f(u(t))-f(u'(t))\bigr)\bigr|
\le \E\,\|u(t)-u'(t)\|_V.
\end{equation}
Recall now that, for any $s\in[1,2]$, we have the interpolation inequality (e.g., see Section~1.6 in~\cite{henry1981})
$$
\|u\|_V\le C_5\|u\|_{L^1}^{1-\theta_s}\|u\|_s^{\theta_s}, \quad u\in H^s,
$$
where we set $\theta_s=\frac{3}{2s+1}$. Combining this with~\eqref{34} and taking the supremum over all $f\in L_b(V)$ with $\|f\|_L\le1$, we derive
$$
\|P_t(v,\cdot)-P_t(v',\cdot)\|_L^*
\le \E\,\Bigl\{\|u(t)-u'(t)\|_{L^1}^{1-\theta_s}\bigl(\|u(t)\|_s+\|u'(t)\|_s\bigr)^{\theta_s}\Bigr\}.
$$
It follows from~\eqref{8} and~\eqref{9} that, for $t\ge1$, the right-hand side of this inequality does not exceed $C_6\|v-v'\|_{L^1}^{1-\theta_s}$, where $C_6>0$ depends only on the $H^1$ norm of the initial functions. It is now straightforward to see that, if~\eqref{33} is satisfied, then~\eqref{6} holds with a sequence~$\delta_m$ going to zero as $m\to\infty$. 

\smallskip
Thus, to complete the proof of stability, it remains to establish~\eqref{9}. In view of Lemma~\ref{l3.1}, it suffices to prove that inequality~\eqref{9} is true for the solution~$v(t,x)$ of problem~\eqref{14}--\eqref{16}. To this end, we repeat essentially the argument of Section~\ref{s3.1}, using the relation (cf.~\eqref{30}) 
\begin{equation} \label{48}
v(t,x)=e^{L}v(t-1)-\int_{t-1}^t e^{(t-r)L}g(r,\cdot)\,\dd r,
\end{equation}
where~$g$ is defined by~\eqref{29}, as well as the inequality
\begin{equation} \label{43}
\|e^{t L}\|_{\LL(H,H^s)}\le C_7t^{-s/2} \quad \mbox{for $t\ge0$},
\end{equation}
where $s\in[1,2]$. Namely, suppose we have proved that
\begin{equation} \label{45}
\E \sup_{t\le r\le t+1}\|v(r)\|_1^2\le C_8
\quad\mbox{for $t\ge0$}.
\end{equation}
Combining this with~\eqref{43}, we obtain the following estimate for the first term in~\eqref{48}:
\begin{equation} \label{49}
\E\,\|e^Lv(t-1)\|_2\le C_9\quad\mbox{for all $t\ge1$}. 
\end{equation}
Denoting by~$G(t)$ the second term in~\eqref{48} and using again~\eqref{43}, we write
\begin{equation} \label{50}
\|G(t)\|_s\le C_7\int_{t-1}^t(t-r)^{-s/2}\|g(r)\|\,\dd r
\le C_{10}\sup_{t-1\le r\le t}\|g(r)\|.
\end{equation}
It follows from the definition of~$g$ that (see~\eqref{29})
$$
\|g(r)\|\le C_{11}\bigl(\|z(r)\|_1^2+\|v(r)\|_1^2\bigr)+\|h\|.
$$
Combining this with~\eqref{50}, \eqref{45}, and~\eqref{18}, we obtain
$\E\,\|G(t)\|_s\le C_{12}$ for $t\ge1$. Recalling~\eqref{48} and~\eqref{49}, we arrive at inequality~\eqref{9} with~$u$ replaced by~$v$. 

\smallskip
Thus, it remain to establish~\eqref{45}. Without loss of generality, we shall assume that $t\ge1$. In view of inequalities~\eqref{31} and~\eqref{17}, we can find $t_0\in[t-1,t]$ such that 
\begin{equation} \label{51}
\E\,\|v(t_0)\|_1^2\le C_{13}. 
\end{equation}
Furthermore, standard argument (see Proposition~2.4.9 in~\cite{KS-book}) combined with inequality~\eqref{17} shows that  
\begin{equation} \label{52}
\E\,\|v(t)\|^k\le C_{14}(k)\quad\mbox{for $t\ge0$, $k\ge1$}. 
\end{equation}
We now take the scalar product in~$L^2$ of~\eqref{14} with~$-2\p_x^2v$. Integrating by parts and using the H\"older inequality, we get
\begin{align}
\p_t\|\p_xv\|^2+2\nu\|\p_x^2v\|^2
&\le 2\bigl((z+v)\p_x(z+v)-h,\p_x^2v\bigr)\label{53} \\
&\le \nu\|\p_x^2v\|^2+2\|\p_xv\|_{L^3}^3+C_{15}R(t),
\notag
\end{align}
where $(\cdot,\cdot)$ denotes the scalar product in~$L^2$, $R(t)=1+\|z(t)\|_\sigma^4$, and $\sigma\in(\frac32,2)$ is a fixed number. 
Using the continuous embedding $L^3\subset H^{1/6}$ and an interpolation inequality for Sobolev spaces, we derive
$$
2\|\p_xv\|_{L^3}^3\le C_{16}\|v\|_{7/6}^3\le C_{17}\|v\|^{5/4}\|\p_x^2v\|^{7/4}
\le \nu\|\p_x^2v\|^2+C_{18}\|v\|^{10},
$$
Substituting this into~\eqref{53} and integrating in time, we get
$$
\sup_{t\le r\le t+1}\|\p_xv(r)\|^2\le \|\p_xv(t_0)\|^2
+C_{18}\int_{t_0}^{t+1}\bigl(\|v(r)\|^{10}+R(r)\bigr)\,dr. 
$$
Taking the mean value and using~\eqref{18}, \eqref{51}, \eqref{52}, we arrive at~\eqref{45}. This completes the proof of the stability property. 

\subsection{Recurrence}
\label{s3.3}
Inequalities~\eqref{10} and~\eqref{11} combined with the independence of the components of the extended process~$(\uuu_t,\IP_\uuu)$ imply that 
\begin{align*}
\IP_\vvv\{\|u_2\|_s\le R,\|u_2'\|_s\le R\} &\ge p^2\quad 
\mbox{for $v,v'\in V$},\\
\IP_\vvv\{(u_{T_m},u_{T_m}')\in\BBB_m\}&\ge \e_m^2\quad
\mbox{for $v,v'\in B_{V\cap H^s}(R)$},
\end{align*}
where $\vvv=(v,v')$. 
It is a standard fact in the theory Markov processes that~\eqref{7} follows from these two inequalities; e.g., see the proof of Proposition~3.3.6 in~\cite{KS-book}. We thus confine ourselves to the proof of~\eqref{10} and~\eqref{11}. 

\subsubsection*{Proof of~\eqref{10}}
The Kolmogorov--Chapman relation implies that it suffices to establish the following two inequalities:
\begin{align}
P_1(v,B_{L^\infty}(K))&\ge p_1 \quad\mbox{for any $v\in V$},\label{41}\\
P_1(v,B_{V\cap H^s}(R))&\ge \tfrac12 \quad\mbox{for any $v\in B_{L^\infty}(K)$},\label{42}
\end{align}
where $s\in(1,2)$ is any fixed number, $p_1>0$ does not depend on~$v$, and $K$ and~$R$ are sufficiently large positive numbers. The proofs of these inequalities are based on a comparison principle and a regularisation property for parabolic equations. 

\smallskip
To prove~\eqref{41}, we write $u=z+v$ and note that, in view of~\eqref{17}, 
$$
\IP\{\|z(1)\|_{L^\infty}\le R\}\to1\quad\mbox{as $R\to\infty$}. 
$$
Therefore it suffices to find $R_1>0$ such that
\begin{equation} \label{71}
\IP\{\|v(1)\|_{L^\infty}\le R_1\}\ge p_2\quad\mbox{for any $u_0\in V$},
\end{equation}
where $p_2>0$ does not depend on~$u_0$, and~$v(t,x)$ stands for the solution of~\eqref{14}--\eqref{16}. Let us fix any $\sigma\in(\frac32,2)$ and, given $\rho>0$, define the event 
$$
\Gamma_\rho=\Bigl\{\sup_{0\le t\le 1}\|z(t)\|_\sigma\le\rho\Bigr\}. 
$$
Since $z$ is a zero-mean Gaussian process in~$V$ with H\"older continuous trajectories, there is a positive function $p(\rho)$ of the variable $\rho>0$ such that $\IP(\Gamma_\rho)\ge p(\rho)$. We claim that, if $\rho>0$ is sufficiently small and~$R_1>0$ is sufficiently large, then
\begin{equation} \label{72}
\|v(1)\|_{L^\infty}\le R_1\quad\mbox{on the event~$\Gamma_\rho$}. 
\end{equation}
Once this property is established, the required inequality~\eqref{71} will follow from the lower bound for $\IP(\Gamma_\rho)$. 

The proof of~\eqref{72} is based on a comparison principle for parabolic equations. We first show that 
\begin{equation} \label{73}
v(1,x)\le R_1\quad\mbox{for all $x\in I$}. 
\end{equation}
To this end, we shall construct a supersolution $v_+(t,x)$ for problem~\eqref{14}--\eqref{16} in the domain $D_1=(0,1)\times I$, that is, a smooth function which is positive  on the lateral boundary of~$D_1$ and greater than $u_0$ for $t=0$, and satisfies the inequality
\begin{equation} \label{74}
\p_tv_+-\nu\p_x^2v_++(v_++z)\p_x(v_++z)-h(x)\ge0
\quad\mbox{for  $(t,x)\in D_1$}. 
\end{equation}
Let us set
$$
v_+(t,x)=\frac{\delta(x+M)+C\e}{t+\e}, 
$$
where $C>\|u_0\|_{L^\infty}$ is fixed, $\e$ and~$\delta$ are small positive parameters, and $M>0$ is a large parameter. It is straightforward to check that 
$v_+(t,0)>0$ and $v_+(t,\pi)>0$ for $t\in[0,1]$ and $v_+(0,x)>\|u_0\|_{L^\infty}$ for $x\in\bar I$. Furthermore, it is a matter of a simple computation to show that, if $\e\in(0,1)$ and $\rho>0$ is sufficiently small, then~\eqref{74} holds, provided that $\delta>0$ is small (for instance, one can take $\delta=\frac12$) and~$M>0$ sufficiently large. Therefore, by Theorem~2.2 in~\cite{AL-1983}, we have 
$$
v(t,x)\le v_+(t,x)=\frac{\delta(x+M)+C\e}{t+\e}\quad\mbox{for $(t,x)\in D_1$}.
$$
Taking $t=1$ and letting $\e\to0^+$, we see that~\eqref{73} holds with $R_1=\delta(\pi+M)$. A similar argument shows that the function $v_-=-v_+$ is a subsolution for~\eqref{14}--\eqref{16} in~$D_1$, and therefore $v(1,x)\ge -R_1$ for $x\in I$.

\smallskip
It remains to prove~\eqref{42}. If the set of initial functions~$u_0$ was bounded in~$B$, inequality~\eqref{9} would imply the required result. However, in our case, it is bounded in~$L^\infty$. To overcome this difficulty, let us remark that, in view of~\eqref{31}, we can find $t_0\in(0,\frac12)$ such that 
$\E\|u(t_0)\|_V^2\le 3C_1$. It follows that, if~$B$ is a ball in~$V$ centred at zero of sufficiently large radius, then 
$$
\IP_{u_0}\{u(t_0)\in B\}\ge 2^{-1/2}\quad\mbox{for any $u_0\in B_{L^\infty}(K)$}. 
$$
The Kolmogorov--Chapman relation and an analogue of inequality~\eqref{9} on the half-line $t\ge\tfrac12$ now implies the required result. 

\subsubsection*{Proof of~\eqref{11}}
Along with~\eqref{0.1}, let us consider the controlled equation
\begin{equation} \label{61}
\p_tu-\nu\p_x^2u+u\p_xu=h(x)+\zeta(t,x), \quad x\in (0,\pi),
\end{equation}
where $\zeta(t)$ is a control function with range in $E:=\lspan\{e_1,e_2\}$. 
Let $\hat u(x)$ be a time-independent solution of problem~\eqref{61}, \eqref{0.2} with $\zeta\equiv0$. Such a solution exists and belong to $H^2\cap V$; e.g., see Section~7 in~\cite[Chapter~I]{lions1969}. The following proposition is a consequence of the main result of~\cite{RS-2015}. 

\begin{proposition} \label{p3.2}
There is $M>0$ such that, given $\e>0$, one can find $T_\e>0$ satisfying the following property: for any $u_0\in V$ there exists  a continuous function $\zeta:[0,T_\e]\to E$ such that
\begin{equation} \label{62}
\|u(T_\e)-\hat u\|_{L^1}<\e, \quad \|u(T_\e)\|_V<M,
\end{equation}
where $u(t,x)$ stands for the solution of~\eqref{61}, \eqref{0.2}, \eqref{0.3}. 
\end{proposition}

A well-known argument based on the above result on approximate controllability (see the proof of Theorem~7.4.1 in~\cite{DZ1996}) enables one to prove that 
\begin{equation} \label{63}
p_{\e,M}(v):=P_{T_\e}\bigl(v,B(\e,M)\bigr)>0\quad\mbox{for any $v\in V$},
\end{equation}
where we set 
$$
B(\e,M)=\{u\in V:\|u-\hat u\|_{L^1}<2\e,\|u\|_V<2M\}. 
$$
Since $B(\e,M)$ is an open subset in~$V$, the Feller property of the Markov process $(u_t,\IP_u)$ in~$V$ implies that the function $p_{\e,M}(v)$ is lower semicontinuous and, hence, separated from zero on compact subsets. In particular, denoting by~$B_m$ the closure of $B(\frac1m,M)$ in~$V$, we can find positive numbers $\e_m>0$ such that
$$
P_{T_m}\bigl(v,B_m\bigr)\ge p_{\frac1m,M}(v)\ge \e_m\quad
\mbox{for $v\in B_{V\cap H^s}(R)$}. 
$$
where $T_m=T_{\frac1m}$. It remains to note that the sets~$B_m$ satisfy~\eqref{33} with $C_4=2M$. This completes the proof of the recurrence property. 

\subsection{Proof of Lemma~\ref{l3.1}}
\label{s3.4}
The fact that almost every trajectory belongs to~$\YY_1$ and inequality~\eqref{17} are well known; see Proposition~2.4.10 and Corollary~2.4.11 in~\cite{KS-book} for the case of the 2D Navier--Stokes system. Thus, we shall prove that $z\in C(\R_+,H^s)$ with probability~$1$ and that~\eqref{18} holds for any $t\ge0$, $s<2$, and $k\ge1$.

To establish the required properties, we first assume that $t=0$. Using the factorisation method described in Section~5.3 of~\cite{DZ1992}, we write~$z$ in the form
\begin{equation} \label{19}
z(t)=\int_0^te^{(t-r)L}y_\alpha(r)(t-r)^{\alpha-1}\dd r,
\end{equation}
where $\alpha\in(0,\frac12)$ is an arbitrary number, $L=\nu\p_x^2$, and
$$
y_\alpha(r)=\sum_{j=1}^2\int_0^re^{(r-\theta)L}e_j(r-\theta)^{-\alpha}\dd\beta_j(\theta). 
$$
As is shown in the proof of Theorem~5.9 in~\cite{DZ1992}, $y_\alpha(r)$ is a Gaussian  process satisfying the inequality
\begin{equation} \label{20}
\E\int_0^T\|y_\alpha(r)\|_1^n\dd r\le C_n(\alpha)T\quad\mbox{for any $T>0$ and $n\ge1$}. 
\end{equation}
Now note that 
\begin{equation} \label{21}
\|e^{tL}\|_{\LL(V,H^s)}\le C t^{(1-s)/2}\quad\mbox{for $s\in[1,2]$ and $t\ge0$}. 
\end{equation}
Taking the $H^s$ norm in~\eqref{19}, using~\eqref{21}, and applying the Young inequality, we obtain
\begin{equation} \label{23}
\|z(t)\|_s\le\int_0^t(t-r)^{\alpha-\frac{1+s}{2}}\|y_\alpha(r)\|_V\dd r
\le \int_0^t\bigl((t-r)^{p_m}+\|y_\alpha(r)\|_V^m\bigr)\,\dd r,
\end{equation}
where $p_m=\frac{m}{m-1}(\alpha-\frac{1+s}{2})$, and the number~$m\ge1$ is chosen below. For a given $s\in(1,2)$, we can find $\alpha$ satisfying the inequality $\frac12>\alpha>\frac{s-1}{2}$. In this case, we have $\alpha-\frac{1+s}{2}>-1$, and one can choose~$m$ so large that $p_m>-1$. Combining this with~\eqref{23} and~\eqref{20}, we arrive at~\eqref{18}. The continuity of~$z$ as a function of time with values in~$H^s$ can now be proven by a simple approximation argument. 

To prove~\eqref{18} in the general case, we can assume that $t=t_0\ge1$. Let us write
\begin{equation} \label{061}
z(t)=e^{(t-t_0-1)L}z(t_0-1)+\sum_{j=1}^2\int_{t_0-1}^te^{(t-r)L}e_j\dd\beta_j(r). 
\end{equation}
The sum on the right-hand side of this relation satisfies the same estimate on the interval $[t_0,t_0+T]$ as the function defined by~\eqref{19} on~$[0,T+1]$. Furthermore, in view of~\eqref{21}, we have the following estimate for the first term on the right-hand side of~\eqref{061}:
\begin{equation} \label{062}
\sup_{t_0\le t\le t_0+T}\|e^{(t-t_0-1)L}z(t_0-1)\|_2\le C\,\|z(t_0-1)\|_V. 
\end{equation}
Recalling~\eqref{17}, we see that the mean value of any degree of the left-hand side of~\eqref{062} is bounded by a constant not depending on~$t_0\ge1$. This completes the proof of the lemma. 

\backmatter
\def\cprime{$'$} \def\cprime{$'$}
  \def\polhk#1{\setbox0=\hbox{#1}{\ooalign{\hidewidth
  \lower1.5ex\hbox{`}\hidewidth\crcr\unhbox0}}}
  \def\polhk#1{\setbox0=\hbox{#1}{\ooalign{\hidewidth
  \lower1.5ex\hbox{`}\hidewidth\crcr\unhbox0}}}
  \def\polhk#1{\setbox0=\hbox{#1}{\ooalign{\hidewidth
  \lower1.5ex\hbox{`}\hidewidth\crcr\unhbox0}}} \def\cprime{$'$}
  \def\polhk#1{\setbox0=\hbox{#1}{\ooalign{\hidewidth
  \lower1.5ex\hbox{`}\hidewidth\crcr\unhbox0}}} \def\cprime{$'$}
  \def\cprime{$'$} \def\cprime{$'$} \def\cprime{$'$}
\providecommand{\bysame}{\leavevmode\hbox to3em{\hrulefill}\thinspace}
\providecommand{\MR}{\relax\ifhmode\unskip\space\fi MR }

\providecommand{\MRhref}[2]{%
  \href{http://www.ams.org/mathscinet-getitem?mr=#1}{#2}
}
\providecommand{\href}[2]{#2}

\end{document}